\documentclass[10pt,a4paper]{article} 
\usepackage{amssymb}

\newcommand{\cd}{D}
\newcommand{\ce}{F}

\newcommand{\E}{\mathbb{E}}
\newcommand{\ee}{\mathrm{e}}
\renewcommand{\P}{\mathbb{P}}

\renewcommand{\le}{\leqslant}
\renewcommand{\ge}{\geqslant}

\newtheorem{theorem}{Theorem}
\newtheorem{definition}{Definition}
\newtheorem{proposition}[definition]{Proposition}

\begin{document}
\title{Asymptotics of iterated branching processes}
\author{Didier Piau}
\date{\em Universit\'e Lyon 1}

\maketitle

\begin{abstract}
We study the iterated Galton-Watson process $(X_n)_n$, possibly with thinning, 
introduced by Gawe\l\ and
Kimmel to model the number of repeats of DNA triplets during some
genetic disorders. If the process involves some thinning, then
extinction $\{X_n\to0\}$ and explosion $\{X_n\to\infty\}$ can have
positive probability simultaneously.  If the underlying (simple)
Galton-Watson process is nondecreasing with mean $m$, then, conditionally on
explosion, $\log X_{n+1}\sim X_n\cdot\log
m$ almost surely. This simplifies arguments of Gawe\l\ and Kimmel, and
confirms and extends a conjecture of Pakes.
\end{abstract}

\bigskip





\section*{Introduction}

The parameters of an iterated Galton-Watson (IGW) process are 
a real number $\vartheta$ in $(0,1]$,
  called the thinning parameter, and a
probability measure  $(p_k)_{k\ge0}$ on the
nonnegative integers, called the reproduction law.
The associated IGW process is a Markov chain
$(X_n)_{n\ge0}$ with nonnegative integer values, whose stochastic
evolution is as follows.

Let $(Z^{(n)})_{n\ge0}$ denote an auxiliary i.i.d.\ collection of
Galton-Watson processes
with common reproducing distribution $(p_k)$ and starting from
$Z^{(n)}_0:=1$. Let $(\varepsilon^{(n)}_k)_{n\ge0,\,k\ge1}$ denote
an auxiliary i.i.d.\ collection 
of Bernoulli random
variables with common distribution
$$
\P(\varepsilon^{(n)}_k=1):=\vartheta=:1-\P(\varepsilon^{(n)}_k=0).
$$
Assume that $X_0$, $(Z^{(n)})$ and $(\varepsilon^{(n)}_k)$ are independent.
Let $S^{(n)}_x$ denote  the
population up to level $x$ of the $n$th Galton-Watson processes
$Z^{(n)}$, that is, 
$$
S^{(n)}_x:=Z^{(n)}_1+\cdots+Z^{(n)}_x.
$$ 
If $X_n=0$, then $X_{n+1}:=0$. If
$X_n=x\ge1$, then
$X_{n+1}:=\varepsilon^{(n)}_1+\cdots+\varepsilon^{(n)}_y$, where
$y:={S^{(n)}_x}$.
In other words, the distribution of $X_{n+1}$ conditionally on
$\{X_n=x\}$ and $S^{(n)}_x$
is binomial with
parameters $S^{(n)}_x$ and $\vartheta$.

Gawe\l\ and Kimmel (1996) introduced IGW processes
to model the explosive growth of the number of repeats of DNA triplets
in specific regions of the genome during heritable
disorders such as fragile X-syndrome.
Here, $1+X_n$ models the length of a linear chain of DNA repeats
after $n$ replications, and what we know of the molecular
mechanism of replication suggests to describe the evolution of $(X_n)$
as above.
(Gawe\l\ and Kimmel use a different convention and
define the IGW process as $(X'_n)_{n\ge0}$ with $X'_n:=1+X_n$.)
Indeed, simulations in Gawe\l\ and Kimmel when $p_0=0$ and
$\vartheta\neq1$ suggest
that IGW  processes either die out or grow extremely fast
 after a period of relative quiescence. This is precisely
the behaviour
of the number of repeats of DNA triplets during some of
these genetic disorders.
See also the book by
Kimmel and Axelrod
(2002), which repeats the analysis of Gawe\l\ 
and Kimmel.

We define the explosion $\ce$ and death $\cd$ of the IGW process
as the events
$$
\ce:=\{X_n\to\infty\},
\quad
\cd:=\{X_n\to0\}=\{X_n=0\ \mbox{for $n$ large enough}\}.
$$
Thus $\ce$ and $\cd$ are mutually exclusive.
(In the context of genetic disorders, the death $\cd$ of the
process corresponds to the extinction of the diseased gene lineage, 
and the explosion $\ce$ to the death of
the patient.)
Gawe\l\ and Kimmel show that $\P(\cd)=1$ as soon as
 $p_0\neq0$, and that $\P(\ce)=1$ when $p_0=0$, $p_1\neq1$ and 
$\vartheta=1$.
Later on, their arguments were simplified
 by Pakes (2003).
Pakes also conjectured that, when $p_0=0$, $p_1\neq1$ 
and $\vartheta=1$,
$(\log X_{n+1})/X_n$ converges almost surely to $\log m$, where
 $m>1$ 
denotes the
mean of $(p_k)$, that is
$$
m:=\sum_{k\ge0}k\,p_k.
$$
In this paper, we determine the asymptotic behaviour of every IGW
process. We confirm the conjecture of Pakes and extend it to
IGW processes with thinning, that is, to the case $\vartheta\neq1$,
and we refine
partial results of Gawe\l\ 
and Kimmel which are not recalled above. 

Propositions~\ref{p.p} and \ref{p.pp} below state 
some simple facts about the mean behaviour of IGW processes and
about their
almost sure behaviour in some degenerate cases.
Some of these are due to Gawe\l\ 
and Kimmel, or to Pakes.
We write $\P_x$ and $\E_x$ for the probability
$\P$ and the expectation $\E$, conditional on $X_0=x$.

\begin{proposition}
\label{p.p}
The mean behaviour of the IGW process is as follows.

{\bf (1)}
If $m>1$, then
$\E_x(X_n)\to\infty$ for every $x\ge1$.

{\bf (2)}
If $m<1$, or if $m=1$ and $\vartheta\neq1$, 
then $\E_x(X_n)\to0$ for every $x\ge0$.

{\bf (3)}
If $m=\vartheta=1$, then $\E_x(X_n)=x$ for every $x\ge0$.
\end{proposition}
\begin{proposition}
\label{p.pp}
The almost sure behaviour of the IGW process 
in some degenerate cases is as follows.

{\bf (4)}
If $p_0\neq0$, then $\P_x(\cd)=1$ for every $x\ge0$.

{\bf (5)}
If $p_0=0$, $p_1\neq1$, and $\vartheta=1$, 
then $\P_x(\ce)=1$ for every $x\ge1$.
\end{proposition}

The hypotheses of statements (1) and (4)
are compatible, 
hence one can have
$\E_x(X_n)\to\infty$ and $\P_x(\cd)=1$ simultaneously,
for every $x\ge1$.
This contrasts with the behaviour of usual Galton-Watson processes.
We turn to the non degenerate case that 
proposition \ref{p.pp} leaves out.

\begin{theorem}
\label{t.a}
Assume that $p_0=0$ and $p_1\neq1$.

{\bf (6)}
For every $x\ge0$,
$\P_x(\cd)+\P_x(\ce)=1$.

{\bf (7)}
Assume furthermore that $\vartheta\neq1$.
For every
$x\ge1$, 
$\P_x(\cd)$ and
$\P_x(\ce)$
are both positive, and
$\P_x(\cd)\le\P_1(\cd)^x$.
Thus, $\P_x(\ce)\to1$ when $x\to\infty$.
\end{theorem}

We stress that (7) does not state that
$\P_x(\cd)$ decreases
geometrically when $x\to\infty$. In fact,
one can show in this case that 
$$
\log1/\P_x(\cd)\gg x,
\quad\mbox{when}\ x\to\infty.
$$
We now state our main result.

\begin{theorem}
\label{t.b}
Assume that $p_0=0$.
Conditionally on the explosion $\ce$,
the random variable $(\log X_{n+1})/X_n$ converges almost surely to $\log m$.
In particular, conditionally on $\ce$, $X_{n+1}/X_n$ converges almost
surely to infinity.  
\end{theorem}

As regards the process conditional on the death $D$, one could try to
show that, if suitably renormalized, the death time, that is, the first
hitting time of $0$, converges in distribution, conditionally on $D$ and
when the starting point of the process goes to infinity. We do not
pursue this in the present paper.

In sections~\ref{s.2} and \ref{s.3}, we (re)prove
propositions~\ref{p.p} and \ref{p.pp}.  In section~\ref{s.4}, we deal
with the easy parts of theorem~\ref{t.a}, that is, (6) and the facts
that the probability of death is not zero and that it is at most
geometric.  In section~\ref{s.5}, we provide explicit upper bounds of
the probability of death in some specific cases.  Section~\ref{s.6}
exposes a strategy of proof for the study of the explosion case. In
section~\ref{s.7}, we apply this strategy, first to the proof of the
remaining assertion of (7), that is, the fact that the probability of
explosion is not zero, and finally to the proof of theorem~\ref{t.b}.

\section{Proof of proposition \ref{p.p}}
\label{s.2}

From the construction of the IGW process, one sees that, for every $x\ge0$,
$\E_x(X_1)=\chi(x)$, where the sequence $\chi$ is defined by
$\chi(0):=0$ and, for every $x\ge1$,
$$
\chi(x):=\vartheta\,(m+\cdots+m^x).
$$
When $m>1$, $\chi$ can be extended to a convex function on
$[0,+\infty)$, which we still call $\chi$. 
(This step of the proof would be false for $m<1$.)
Thus,
for every $x\ge0$ and every $n\ge0$,
$$
\E_x(X_{n+1})\ge\chi(\E_x(X_n)).
$$
Choose $x_0\ge1$ large enough, such that $\chi(x_0)\ge 2x_0$.
Then, $\E_{x_0}(X_n)\ge x_0\,2^n$ for every $n\ge0$. (This step of the
proof
uses the
fact that $\chi$ is nondecreasing).
For smaller values of $x$, for instance for $x=1$, consider the event
that $X_n\ge x_0$. For $n$ large enough, $n=n_0$ say, 
this event has positive probability $v_0$
with respect to $\P_1$.
Finally, (1) holds since,
for every $n\ge n_0$ and every $x\ge1$,
$$
\E_x(X_{n})\ge\E_1(X_{n})\ge\P(X_{n_0}\ge x_0)\,\E_{x_0}(X_{n-n_0})
\ge v_0\,x_0\,2^{n-n_0}.
$$
When $m=1$, $\E_x(X_1)=\vartheta\,x$ thus the $m=1$ cases in
 (2) and (3) are obvious.
Finally, when $m<1$,
$\E_x(X_1)\le\vartheta\,m\,x$ and $\vartheta\,m<1$, thus
$\E_x(X_n)\to0$.
This completes the proof of (2).

\section{Proof of proposition \ref{p.pp}}
\label{s.3}

The proof of (5) is direct since in that case, $X_{n+1}\ge X_n$ almost surely.

As regards (4),
if $p_0\neq0$, for each $n\ge0$, the event that $Z^{(n)}_1=0$ has positive
probability $p_0$ and these events
are independent.
Almost surely,
one of them is realized, say $Z^{(n)}_1=0$.
Then $S^{(n)}_x=0$ for every $x\ge1$ and $X_{n+1}=0$, which proves (4).
The argument shows also that 
$$
\P_x(X_n\neq0)\le(1-p_0)^n,
$$ 
for
every $n$ and $x$.
In other words, the time to absorption is 
at most geometric, uniformly.

\section{Elementary parts of theorem \ref{t.a}}
\label{s.4}
\subsection{Almost surely, death or explosion}

This is straightforward, and analogous to the usual Galton-Watson case.
One has to check that $0$ is the only non transient state of the
Markov chain $(X_n)$.
If $p_0=0$ and $\vartheta=1$, this is true, see the proof of (2).
If $p_0=0$ and $\vartheta\neq1$, the return to $x\ge1$, starting from
$x$,  assumes that the first step is not to $0$.
Thus, it has a probability at most $1-\P_x(X_1=0)<1$, see section \ref{ss.2}.

\subsection{Probability of death, not zero}
\label{ss.2}

When 
$X_0=x$, $X_1=0$ iff the
thinning kills each and every $S^{(0)}_x$ individuals.
Hence, for
every $x\ge0$,
$$
\P_x(\cd)\ge\P_x(X_1=0)=\E((1-\vartheta)^{S_x}),
$$
which is positive.
\subsection{Probability of death, at most geometric}
\label{s.43}

Assume that $p_0=0$ and fix $x\ge1$ and $y\ge1$.
The fundamental
branching property of the Galton-Watson process $Z^{(n)}$
means that
$Z^{(n)}_{x+y}$ is the sum of $Z^{(n)}_x\ge1$ random variables
distributed like $Z^{(n)}_y$. Hence $S^{(n)}_{x+y}$ is stochastically
greater than the sum of $S^{(n)}_{x}$ and of an independent copy of
$S^{(n)}_{y}$.
After thinning, this shows that $X_1$ under $\P_{x+y}$ is stochastically
greater than the sum of $X_1$ under $\P_x$ and of an independent copy of
$X_1$ under $\P_y$.
By recursion over $n\ge1$, the same assertion holds when one replaces
$X_1$ by $X_n$.
Hence,
$$
\P_{x+y}(X_n=0)\le\P_{x}(X_n=0)\,\P_{y}(X_n=0).
$$
This implies that $\P_x(\cd)\le\P_1(\cd)^x$ for every
$x\ge0$, an assertion of (7).

Note that the important step here is to prove that $\P_1(\cd)\neq1$.
We do this in section \ref{s.5} in some specific cases, the general
case is in section \ref{s.71}.

\section{Upper bounds of the probability of death}
\label{s.5}

When $p_0=0$, in the special case $m\,\vartheta>1$,
one can bound explicitly $\P_1(\cd)$
by some $q<1$, using the generating function of $(p_k)$.
Thanks to section \ref{s.43}, this proves that 
$\P_x(\cd)\le q^x$. 

\subsection{Special case}
\label{s.51}

One can often get an upper bound of $q:=\P_1(\cd)$ at small cost, as follows.
Since the sequence $(\P_x(\cd))_{x\ge0}$ is submultiplicative,
$$
q=\E_1(\P_{X_1}(\cd))\le\E_1(q^{X_1}).
$$
The generating function of $X_1$ is
$g(s):=\E_1(s^{X_1})=f(1-\vartheta+\vartheta\,s)$, with
$$
f(s):=\sum_{k\ge1}p_k\,s^k.
$$
Thus, $q\le g(q)$ and it is not hard to see that $q$ is at most the
smallest root of the equation $s=g(s)$. 
In the supercritical case $g'(1^-)=m\,\vartheta>1$,
$s\le g(s)$ for $s=1$ and for $s\le q_*$, 
where $q_*$ in $(0,1)$ is defined by 
$q_*=g(q_*)$. Finally, $q\le q_*$.
In particular, $\P_x(\cd)\neq1$.

\subsection{Binary case}
\label{s.52}

Assume that the Galton-Watson process describes a binary
replication with efficiency $\lambda$.
Thus, $f(s):=(1-\lambda)\,s+\lambda\,s^2$ and $m=1+\lambda$.
Elementary computations then yield the following.
If $\vartheta>1/m$, $q\le q(\lambda,\vartheta)$ with
$$
q(\lambda,\vartheta):=
(1-\vartheta)(1-\lambda\,\vartheta)/(\lambda\,\vartheta^2).
$$
Note that $q(\lambda,\vartheta)$ is in $(0,1)$, except when 
$\vartheta=1$, and then $q=q(\lambda,1)=0$, and when 
$\vartheta=1/m$, and then $q(\lambda,\vartheta)=1$ but $q<1$.

\section{A general strategy}
\label{s.6}

Assume that $p_0=0$.  We first explain our strategy for the study of
the explosion of the IGW process.  We fix an integer valued sequence
$(\varphi(x))_{x}$ with $\varphi(x)\ge x$ for every $x$, and an
integer valued sequence $(\psi(x))_{x}$.

\begin{definition}
Let $C$ denote the
event that $X_{n+1}\ge\varphi(X_n)$ for every $n\ge0$.
For every $x\ge1$, let
$\eta_x:=\varepsilon_1+\cdots+\varepsilon_{\psi(x)}$ denote a random
variable of binomial distribution of parameters $\psi(x)$ and
$\vartheta$.
Finally, introduce the probabilities
$$
A(x):=\P(S_x\le\psi(x)),
\quad
B(x):=\P(\eta_{x}\le\varphi(x)).
$$
\end{definition}

To show that $\P_x(C)\neq0$,
we start from
$$
\P_x(X_1\le\varphi(x))\le A(x)+B(x).
$$
From  Markov
inequality and from the fact that $S_x\ge Z_x$,
$$
A(x)\le\P(Z_x\le\psi(x))\le\psi(x)\,\E(1/Z_x)
\le\psi(x)\,\ee^{-c\,x},
$$
with a positive $c$ that depends only on the distribution $(p_k)$, see
for instance proposition A.2 in the appendix of Piau (2004).
Likewise, Chebychev inequality for Bernoulli random variables yields
$$
B(x)\le\ee^{\varphi(x)}\,\E(\ee^{-\varepsilon})^{\psi(x)}
\le\ee^{\varphi(x)-c'\,\psi(x)},
$$
with a positive $c'$ that depends only on $\vartheta$ 
(for instance $c':=\vartheta^2$).
Assume that there exists a nonincreasing sequence  $(\gamma(x))_{x}$,
such that
$$
\P_x(X_1\le\varphi(x))\le\gamma(x)<1,
$$
for every $x$, for instance because
$$
\psi(x)\,\ee^{-c\,x}+\ee^{\varphi(x)-c'\,\psi(x)}\le\gamma(x),
$$
and define recursively $(\gamma_k(x))_{k\ge0}$ by 
$\gamma_0(x):=\gamma(x)$ and 
$\gamma_{k+1}(x):=\gamma_k(\varphi(x))$.
Conditioning successively on the values of $X_n$ and iterating the
above yields 
$$
\P_x(C)\ge\prod_{k\ge0}(1-\gamma_k(x)).
$$
The iteration uses both the fact that $\varphi(x)\ge x$ and the fact that
$(\gamma(x))_x$ is nonincreasing.
As  a consequence, 
when $(\gamma_k(x))_{k}$ is summable and when every $\gamma_k(x)<1$,
the infinite product
is positive. Finally, since $C\subset\ce$,
$\P_x(\ce)\ge\P_x(C)$ and $\P_x(\ce)$ is not zero.

Using standard zero-one laws, this proves at the same time that,
conditionally on $\ce$, the event that $\liminf X_n/n\ge1$ is almost sure,
since this is an asymptotic event which 
contains $C$.
We refine this below.

\section{Some applications of the strategy}
\label{s.7}

\subsection{Probability of explosion, not zero}
\label{s.71}

We first apply section \ref{s.6}
with $\varphi(x):=x+1$ and $\psi(x):=x^2$.
One can choose $\gamma(x)<1$ for every $x\ge1$ and such that
$\gamma(x)\le\ee^{-c''\,x}$ when $x\to\infty$, with a positive $c''$.
Since $\gamma_k(x)=\gamma(k+x)$,
the series  $(\gamma_k(x))_{k}$ is summable and  every 
$\gamma_k(x)<1$.
This proves that $\P_x(\ce)$ is positive for every $x\ge1$, that is, the
missing part  of
(7), and completes the proof of theorem~\ref{t.a}.

\subsection{Lower bound in theorem \ref{t.b}}
\label{s.72}

Our second application of section \ref{s.6} is more involved.
We choose  $\mu<m$ and an integer sequence $(\varphi(x))_{x}$ such
that
$\varphi(x)\sim\mu^x$ for $x$ large enough.
Then we choose $\nu$ in $(\mu,m)$ and an integer sequence $(\psi(x))_x$
such that
$\psi(x)\sim\nu^x$  for $x$ large enough.
Since $\nu>\mu$, the contributions $(B(x))_x$ are summable.

As regards the contributions $(A(x))_x$, we use standard
estimates of the harmonic moments of $Z_x$.
Choose a positive $r$ such that $p_1\,m^r<1$, this is possible as soon
as $m$ is finite and $p_1\neq1$.
From Ney and Vidyashankar (2003),
the behaviour of
$\E(1/(Z_x)^r)$ is ruled by the so-called Seneta constants. 

To keep
things simple, we choose $\rho$ in $(\nu,m)$ and we use the following 
easy consequence
of the results by Ney and Vidyashankar. There exists a finite constant
$c_0$ such that $\E(1/(Z_x)^r)\le c_0\,\rho^{-rx}$ for every $x\ge1$.
Then, the Markov
inequality for $(Z_x)^r$ yields, for $x$ large enough,
$$
A(x)\le\P(Z_x\le\psi(x))\le\psi(x)^r\,\E(1/(Z_x)^r)\le c_0\,(\nu/\rho)^{rx}.
$$
Finally, the sequence $(\gamma(x))_{x}$ is allowed to decrease
geometrically, hence to be
summable.

Since $\varphi(x)\ge x+1$ for $x$ large enough, 
$\gamma_k(x)\le\gamma(x+k)$ and the
  series $(\gamma_k(x))_{k}$ is summable for every $x$.
Thus $\P_x(C)\neq0$. 
Introduce
$$
G:=\{\liminf Y_{n}\ge\log\mu\},
\quad
Y_{n}:=(\log X_{n+1})/X_n.
$$
Since $C\subset G$ and $G$ is asymptotic,
$G$ is almost sure on $\ce$.
Finally, conditionally on $\ce$,
$$
\liminf Y_{n}\ge\log m
\quad
\mbox{almost surely}.
$$

\subsection{Upper bound in theorem \ref{t.b}}
\label{s.73}

As regards the other side of the equality of theorem \ref{t.b}, fix $\mu>m$. We
use the simple fact that, if $Y_n\ge\log\mu$, then, conditionally on $X_n=x$,
one has $S_x^{(n)}\ge X_{n+1}\ge\mu^x$.
Furthermore,
$$
\P(S_x\ge\mu^x)\le\mu^{-x}\,\E(S_x)\le(m/\mu)^x\,m/(m-1),
$$
which is summable.
Conditionally on $\ce$, $\liminf X_n/n\ge1$ is almost sure, see the last
lines of section \ref{s.6}.
This implies that
the events that $Y_n\ge\log\mu$
are realized at most for a finite number of values of $n$, 
conditionally on $\ce$ or
not. Thus, 
$\limsup Y_{n}\le\log\mu$ almost surely.
This concludes the proof of theorem~\ref{t.b}.

\bigskip

Universit\'e Claude Bernard Lyon 1 \\ Institut Camille Jordan UMR 5208
\\ 
Domaine de Gerland \\
50, avenue Tony-Garnier \\ 69366 Lyon Cedex 07 (France)

{\tt
Didier.Piau@univ-lyon1.fr} \\ {\tt http://lapcs.univ-lyon1.fr}

\end{document}